\documentclass{article}
\usepackage[french]{babel}
\usepackage{abou}
\usepackage{graphicx}

\input tcilatex
\begin{document}

\label{toto}

\begin{center}
{\LARGE {Forme des connexes de Farey} \\[0pt]
}{\Large {Par Saab Abou-Jaoud\'{e}\footnote{{\Large ancien professeur de
Math\'{e}matiques Sp\'{e}ciales}}}}
\end{center}

{\Large \ }

{\Large \hskip-1cm \textbf{R\'{e}sum\'{e} : } }

L'objectif de ce texte est de d\'{e}montrer que les composantes connexes du
complexe de Farey plan sont des triangles ou des quadrilat\`{e}res. Pour le
faire nous faisons un retour sur les polygones convexes plan \`{a} bord
orient\'{e} en d\'{e}montrant que si deux vecteurs c\^{o}t\'{e}s
cons\'{e}cutifs du bord d'un tel polygone ne sont jamais dans le m\^{e}me
quadrant, alors ce polygone est un triangle ou un quadrilat\`{e}re. Ce
r\'{e}sultat est, \`{a} notre connaissance, in\'{e}dit. Nous appliquons ce
r\'{e}sultat aux connexes de Farey en d\'{e}montrant que les polyg\^{o}nes
qui les d\'{e}limitent la v\'{e}rifient.

{\Large \hskip-1cm \textbf{Mots cl\'{e} :} }

Complexe de Farey, Connexe de Farey, Polygone convexe, Convexe polygonal
direct, Droite orient\'{e}e, Demi-plan.

\section{Introduction}

La conjecture de Tajine-Daurat trouve sa source en g\'{e}om\'{e}trie
discr\`{e}te au laboratoire du professeur Tajine \`{a} Strasbourg. Elle
s'\'{e}nonce de la mani\`{e}re suivante :

Soit $(m,n)\in \N^{*2}.$ On consid\`{e}re l'ensemble $D_{m,n}$ des droites
rencontrant le carr\'{e} unit\'{e} $CU=[0,1]^{2}$ et d'\'{e}quation$%
\;ux+vy-w=0$\ dont les coefficients $u,v,w$ sont des \underline{entiers
relatifs} et v\'{e}rifient les conditions :$\;\;$%
\[
|u|\le m,\;|v|\le n,\;(u,v)\ne (0,0),\;w\in \Z
. 
\]
On consid\`{e}re l'ensemble $A,$ compl\'{e}mentaire dans le plan de l'union
des droites de $D_{m,n}.$ Alors les composantes connexes born\'{e}es de $A$
contenues dans le carr\'{e} unit\'{e} $CU$ sont des triangles ou des
quadrilat\`{e}res. Nous donnons une condition pour qu'une telle composante
connexe $K$ soit un triangle. En particulier si l'un des sommets de $K$ a
pour coordonn\'{e}e $(p/q,p^{\prime }/q^{\prime }),$ avec $0<q\le m$ ou$%
\;0<q^{\prime }\le n,$ alors $K$ est un triangle.

Pour faire cette \'{e}tude nous aurons besoin de faire un long d\'{e}tour
passant par les polygones convexes et les convexes polygonaux direct, car
nous nous sommes aper\c{c}u que la propri\'{e}t\'{e} ci-dessus est un cas
particulier d'un th\'{e}or\`{e}me plus g\'{e}n\'{e}ral concernant les
polygones convexes et qui, \`{a} notre connaissance, est \`{a} ce jour
in\'{e}dit.

Nous commen\c{c}ons par rappeler, dans la section II, les d\'{e}finitions
des objets g\'{e}om\'{e}triques utilis\'{e}s. Dans la section III, nous
donnons les d\'{e}finitions des notions de polygone convexe et de convexe
polygonal direct et nous rappelons le th\'{e}or\`{e}me classique de
dualit\'{e} entre ces deux notions. Nous introduisons \'{e}galement la
notion de r\'{e}duction d'un convexe polygonal direct. Dans la section IV,
nous d\'{e}montrons le th\'{e}or\`{e}me principal concernant les polygones
convexes. Nous d\'{e}finissons dans la section V l'ensemble des droites $%
D_{m,n}$ et le complexe de Farey $CF(m,n)$ et nous d\'{e}montrons le
th\'{e}or\`{e}me de structure des composantes connexes de $CF(m,n)$. Enfin,
dans la section VI, nous \'enon\c{c}ons la conjecture forte de Tajine-Daurat.

\section{Rappel des d\'{e}finitions utiles.}

On est dans le plan affine r\'{e}el orient\'{e} $E$, muni d'un rep\`{e}re $%
(O,i,j)$\footnote{{On notera les vecteurs sans fl\^{e}che au dessus pour des
raisons de commodit\'{e}s, sauf n\'{e}cessit\'{e}.}} direct. On le suppose
muni de sa topologie usuelle. La donn\'{e}e du rep\`{e}re $(O,i,j)$ nous
permet de rep\'{e}rer un point $M$ du plan par les coordonn\'{e}es $(x,y)$
du vecteur $OM=xi+yj$. Les coordonn\'{e}es d'un point $M$ seront not\'{e}s $%
(x_{M},y_{M}),$ lorsque plusieurs points sont en jeu. On note $\det (V,W)$
le d\'{e}terminant dans la base $(i,j)$ des vecteurs $V=ai+bj$ et $W=ci+dj.$
On a donc : 
\[
\det (V,W)=\left| 
\begin{array}{cc}
a & c \\ 
b & d
\end{array}
\right| =ad-bc. 
\]

Les droites $(O,i)$ et $(O,j)$ d\'{e}finissent dans le plan priv\'{e} de $O$
quatre quadrants : 
\begin{eqnarray*}
Q_{1} &=&\{(x,y)\;|\;x\ge 0,\;y\ge 0\},\;Q_{2}=\{(x,y)\;|\;x\le 0,\;y\ge 0\},
\\
Q_{3} &=&\{(x,y)\;|\;x\le 0,\;y\le 0\},\;Q_{4}=\{(x,y)\;|\;x\ge 0,\;y\le 0\}.
\end{eqnarray*}
On peut ainsi dire si un point $M$ appartient \`{a} l'un ou l'autre de ces
quadrants (d'origine $O$) en disant que le vecteur $OM$ lui appartient. Si $%
A $ est un point du plan, on d\'{e}finit de mani\`{e}re analogue les
quadrants d'origine $A.$ Notons d\'{e}j\`{a} qu'un point peut appartenir
\`{a} deux quadrants s'il est situ\'{e} sur l'un des axes de coordonn\'{e}es
et qu'il est distinct de $O$. On se permettra d'\'{e}tendre cette notion aux
vecteurs $V$ du plan vectoriel.

Nous aurons besoin de la notion de droite orient\'{e}e d\'{e}crite par un
couple $(A,B)$ de points. La notation $(A,B)$ d\'{e}signera l'unique droite
du plan contenant les points $A$ et $B$ lorsque ces deux points sont
distincts. La droite $(A,B)$ admet des \'{e}quations de la forme $ux+vy-w=0$%
, mais elles sont toutes proportionnelles. Parmi ces \'{e}quations nous en
distinguons une que nous qualifierons de \underline{canonique}. Son premier
membre sera : 
\[
f(M)=\det (AB,AM)=(x_{B}-x_{A})(y_{M}-y_{A})-(y_{B}-y_{A})(y_{M}-y_{A}), 
\]
de sorte que : 
\[
M\in (A,B)\Leftrightarrow f(M)=0. 
\]
La fonction $f$ est clairement une fonction affine de $M$.

Cette d\'{e}finition canonique nous permet de d\'{e}finir des r\'{e}gions
dans le plan \`{a} l'aide du signe de $f(M)$. Si $R\in \{\le ,\ge ,<,>\}$,
on notera de mani\`{e}re g\'{e}n\'{e}rale $(A,B)^{R0}$ l'ensemble $%
\{M\;|\;f(M)\,R\,0\}$. Ainsi $(A,B)^{>0}$ (resp. $(A,B)^{<0},$ resp. $%
(A,B)^{\ge 0},$ resp. $(A,B)^{\le 0}$) sera le c\^{o}t\'{e} strictement
positif (resp. strictement n\'{e}gatif, resp. positif, resp. n\'{e}gatif) de
la droite $(A,B)$.

Les trois ensembles $(A,B),\;(A,B)^{>0},\;(A,B)^{<0}$ forment de fa\c{c}on
\'{e}vidente une partition du plan en trois parties toutes convexes, la
premi\`{e}re \'{e}tant ferm\'{e}e, les deux autres \'{e}tant ouvertes. La
fronti\`{e}re de $(A,B)^{>0}$ ainsi que celle de $(A,B)^{\ge 0},\;(A,B)^{<0}$
et $(A,B)^{\le 0}$ est la droite $(A,B)$. Notons enfin que $%
(B,A)^{>0}=(A,B)^{<0}$ (de m\^{e}me $(B,A)^{\ge 0}=(A,B)^{\le 0}$).

Nous aurons besoin de la notion de segment associ\'{e} \`{a} une couple de
points $(A,B)$ qu'on notera $[A,B]$. Cet ensemble est d\'{e}fini par la
formule : 
\[
\lbrack A,B]=\{(1-t)A+tB\;|\;t\in [0,1]\}. 
\]
Quand $A=B$, Le segment $[A,B]$ est r\'{e}duit \`{a} un point. Dans la suite
nous ferons l'abus d'appeler ``segment'' un segment non r\'{e}duit \`{a} un
point. Cela nous \'{e}vitera des lourdeurs.

Voici une remarque importante concernant ces deux notions. Si la droite $%
(A,B)$ traverse le segment $[C,D]$ (i.e. le coupe en un point $M\notin
\{C,D\}$), alors $C$ appartient \`{a} l'un des demi-plans ouvert $(A,B)^{>0}$
ou $(A,B)^{<0}$ et $D\ $appartient \`{a} l'autre. On dira alors que $C$ et $%
D $ sont de part et d'autre de la droite $(A,B)$.

R\'{e}ciproquement, si $f(M)=0$ est une \'{e}quation de la droite $(A,B)$ et
si $C$ et $D$ sont deux points du plan qui v\'{e}rifient $f(C)>0,\;f(D)\le
0, $ alors, par le th\'{e}or\`{e}me des valeurs interm\'{e}diaires
appliqu\'{e} \`{a} la fonction continue $t\rightarrow f\left(
(1-t)C+tD)\right) ,$ la droite $(A,B)$ rencontre le segment $[C,D]$ en un
point diff\'{e}rent de $C$. Si la condition $f(D)\le 0$ est remplac\'{e}e
par la condition $f(D)<0$, alors la droite $(A,B)$ traverse le segment $%
[C,D] $. Cette remarque nous servira dans la d\'{e}monstration du
th\'{e}or\`{e}me \ref{t1}\\[0pt]

\section{Polygone convexe et convexe polygonal direct.}

\subsection{D\'{e}finition d'un polygone convexe ($PC$).}

On dit qu'un partie $K\subset E$ est un $PC$ (Polygone Convexe) s'il existe
un ensemble fini $X$ de demi-plans ferm\'{e}s tels que :

\begin{enumerate}
\item  {\ $K=\bigcap\limits_{P\in X}P.$ }

\item  {\ $K$ est born\'{e}. }

\item  {\ L'int\'{e}rieur de $K$ est non vide. }
\end{enumerate}

On dira que $X$ est un g\'{e}n\'{e}rateur de $K$.

On note $Fr(K)$ la fronti\`{e}re de $K$. Elle est contenue dans l'union des
droites fronti\`{e}res des $P\in X$. Parmi les g\'{e}n\'{e}rateurs de $K$,
il y en a un $X_{0}$ qui est de cardinal minimum. Ce cardinal sera le nombre
de c\^{o}t\'{e}s de $K$. On le note $NC(K)$ (Nombre de C\^{o}t\'{e}s). On
remarquera (et cela n'est pas \'{e}vident) que $P\in X_{0}$ implique que $%
K\cap Fr(P)$ est un segment (sinon cela contredit la minimalit\'{e} du
cardinal de $X_{0},$ la suppression de $P$ de l'ensemble $X_{0}$ ne
modifiant pas $K$).

Soit $P$ un \'{e}l\'{e}ment de $X$. La fronti\`{e}re de $P$ est une droite
orient\'{e}e $D$. On peut supposer, quite \`{a} changer l'orientation de $D$
par changement du signe de son \'{e}quation, que $P=D^{\ge 0}$. On peut
ainsi remplacer la donn\'{e}e d'un ensemble fini de demi-plans ferm\'{e}s $X$
par la donn\'{e}e d'un ensemble fini $X^{\prime }$ de droites orient\'{e}s
telles que : 
\[
K=\bigcap\limits_{D\in X^{\prime }}D^{\ge 0}. 
\]

Un vrai triangle $T$ est l'exemple type d'un polygone convexe et $NC(T)=3$.

\subsection{D\'{e}finition d'un convexe polygonal direct \`{a} $n$
c\^{o}t\'{e}s ($CPD$).}

Soit $n$ un entier naturel sup\'{e}rieur ou \'{e}gal \`{a} $3$ et $\mathcal{A%
}=(A_{k})_{k\in \Z}$ une famille de points de $E$ index\'{e}e par $\Z$. On
dit que $\mathcal{A}$ est un $CPD$ (Convexe Polygonal Direct) \`{a} $n$
sommets si les conditions suivantes sont r\'{e}alis\'{e}es :

\begin{enumerate}
\item  {\ La suite $\mathcal{A}$ est $n$-p\'{e}riodique (i.e. $\forall p\in %
\Z,\;A_{p+n}=A_{p}$) }

\item  {\ $\forall p\in [0..n-1],$ $\forall q\in [1..n-1]\backslash
\{p,p+1\},\;\det (\overrightarrow{A_{p}A_{p+1}},\overrightarrow{A_{p}A_{q}}%
)>0.$ }
\end{enumerate}

Cette derni\`{e}re condition implique que, $\forall (p,q)\in
[1..n-1]^{2},\;p\ne q\Rightarrow A_{p}\ne A_{q}.$ L'entier naturel $n$ est
donc la p\'{e}riode principale de $\mathcal{A}$. Il suffira de se donner les
points de $\mathcal{A}$ d'indice $k$ tel que $0\le k<n$ pour d\'{e}crire le $%
CPD$. Nous utiliserons cette libert\'{e} en parlant du $CPD$ $(A_{k})_{0\le
k\le n-1}$,

A un $CPD$ $\mathcal{A}=(A_{k})_{0\le k\le n-1}$, on associe $LP(\mathcal{A}%
)=\bigcup\limits_{0\le k\le n-1}[A_{k},A_{k+1}]$ qui est la ligne polygonale
associ\'{e}e \`{a} $\mathcal{A}$. On lui associe \'{e}galement la suite des
droites orient\'{e}es $\left( (A_{k},A_{k+1})\right) _{0\le k\le n-1}$.

Les trois sommets d'un vrai triangle num\'{e}rot\'{e}s dans le sens
convenable est l'exemple type d'un $CPD$.

\subsection{Relation entre $PC$ (Polygone Convexe) et $CPD$ (Convexe
polygonal Direct).}

Entre ces deux notions, il existe une relation bien connue, \'{e}nonc\'{e}e
dans le th\'{e}or\`{e}me suivant qui permet de remplacer l'une des notions
par l'autre :

\begin{Th}
- {Soit $K$ un $PC$. Alors il existe un $CPD$\ not\'{e} $\mathcal{A}$ tel
que $Fr(K)=LP(\mathcal{A})$. }\newline
- {Soit $\mathcal{A}$ un $CPD$. Alors il existe un $PC$ not\'{e} $K$ tel que 
$LP(\mathcal{A})=Fr(K)$. }
\end{Th}

Ce th\'{e}or\`{e}me est bien connu des g\'{e}om\`{e}tres. En l'absence de
r\'{e}f\'{e}rence, nous en avons fait une d\'{e}monstration fond\'{e}e sur
l'existence d'un g\'{e}n\'{e}rateur minimal de $K$, mais Monsieur Fran\c{c}%
ois Rideau nous a signal\'{e} la pr\'{e}sence d'une formulation
\'{e}quivalente dans le c\'{e}l\'{e}bre livre de g\'{e}om\'{e}trie de Marcel
Berger [1], un grand classique.

Pourquoi tout cela ? parce que nous aurons besoin de r\'{e}duire un $CPD$
par suppression d'un point et que ce qui reste soit encore un $CPD$. Comme
c'est un point cl\'{e} de notre raisonnement, nous en donnons une
d\'{e}monstration dans le paragraphe qui suit.

\subsection{R\'{e}duction d'un $CPD$.}

Soit $\mathcal{A}=(A_{k})_{k\in \Z}$ un $CPD$ \`{a} $n$ sommets, avec $n\ge
4.$ On d\'{e}finit une suite $(B_{k})_{0\le k\le n-2}$ Par la formule $%
B_{k}=A_{k+1}$ avec $0\le k\le n-2.$ Cette suite $(B_{k})$ est prolong\'{e}e
\`{a} $\Z$ par $(n-1)$-p\'{e}riodicit\'{e}.

\begin{Th}
{\label{t1}La suite $(B_{k})_{k\in \Z}$ d\'{e}finie ci-dessus est un $CPD$ }
\end{Th}

\begin{Dem}
Notons $K$ le polygone convexe associ\'{e} \`{a} $\mathcal{A}$. Pour
d\'{e}montrer ce th\'{e}or\`{e}me, il nous suffit d'\'{e}tablir que tous les
sommets $(A_{k})_{2\le k\le n-2}$ sont dans le demi-plan ouvert $%
P=(A_{n-1},A_{1})^{>0}$. Notons $A=A_{0},$ $B=A_{1},$ $C=A_{n-1}$ et pla\c{c}%
ons-nous dans le rep\`{e}re direct $(A,\overrightarrow{AB},\overrightarrow{AC%
})$. Ce rep\`{e}re est bien direct car : 
\[
\det (\overrightarrow{AB},\overrightarrow{AC})=\det (\overrightarrow{%
A_{0}A_{1}},\overrightarrow{A_{0}A_{n-1}})>0. 
\]
Il est facile de voir que $P=(A_{n-1},A_{1})^{>0}=\{\;(x,y)\;|\;x+y-1>0\;\}%
\footnote{{En effet, $\det (\overrightarrow{A_{n-1}A_{1}},\overrightarrow{%
A_{n-1}M})=\left| 
\begin{array}{cc}
1 & x \\ 
-1 & y-1
\end{array}
\right| =y-1+x.$}}.$\newline
Notons $(a_{k},b_{k})$ les coordonn\'{e}es du point $A_{k},\;2\le k\le n-2.$
Les hypoth\`{e}ses $\det (A_{0}A_{1},A_{0}A_{k})=\left| 
\begin{array}{cc}
1 & a_{k} \\ 
0 & b_{k}
\end{array}
\right| >0$ et $\det (\overrightarrow{A_{n-1}A_{0}},\overrightarrow{%
A_{n-1}A_{k}})=\left| 
\begin{array}{cc}
0 & a_{k} \\ 
-1 & b_{k}-1
\end{array}
\right| >0$ montrent que $b_{k}>0$ et $a_{k}>0.$L'hypoth\`{e}se $\det (%
\overrightarrow{A_{1}A_{2}},\overrightarrow{A_{1}A_{n-1}})>0$ montre que : 
\[
\left| 
\begin{array}{cc}
a_{2}-1 & -1 \\ 
b_{2} & 1
\end{array}
\right| =a_{2}-1+b_{2}>0, 
\]
ce qui montre que $A_{2}\in P.$ Supposons l'existence d'un entier $k\in
[2,n-2]$ tel que $a_{k}+b_{k}-1\le 0$ et soit $p$ le plus petit d'entre eux.
On a $p\ge 3$ et $p-1$ est tel que $a_{p-1}+b_{p-1}-1>0.$ Puisque, par
d\'{e}finition de l'entier $p,\;a_{p}+b_{p}-1\le 0,$ le segment $%
[A_{p-1},A_{p}]$ rencontre la droite $(B,C),$ dont l'\'{e}quation est $%
f(x,y)=x+y-1=0,$ en un point $M$ de coordonn\'{e}es $(x_{M},y_{M})$
strictement positives puisque celles de $A_{p-1}$ et celles de $A_{p}$ le
sont. De plus, comme $x_{M}+y_{M}=1,\;$on a $M\ne A_{p-1},$ donc le vecteur $%
\overrightarrow{A_{p-1}M}$ est un vecteur directeur de la droite
orient\'{e}e $(A_{p-1},A_{p}).$ Le point $M,$ \'{e}tant un point de la
droite $(B,C)$ \`{a} coordonn\'{e}es strictement positives, est un point du
segment $[B,C]$ distinct de $B$ et de $C$. Dans ces conditions les deux
d\'{e}terminants $\det (\overrightarrow{A_{p-1}M},\overrightarrow{MB})$ et $%
\det ($ $\overrightarrow{A_{p-1}M},\overrightarrow{MC})$ sont non nuls et de
signe contraire. Les points $B$ et $C$ sont ainsi de part et d'autre du
c\^{o}t\'{e} $(A_{p-1},A_{p})$ de $\mathcal{A},$ ce qui contredit la
convexit\'{e} de $K$.
\end{Dem}

Nous avons accessoirement d\'{e}montr\'{e} ce que nous appelons le lemme des
diagonales :

\begin{Lem}[des diagonales]
Pour $2\le k\le n-2,$ les segments $[A_{0},A_{k}]$ et $[A_{n-1},A_{1}]$ se
traversent.
\end{Lem}

Cette op\'{e}ration de r\'{e}duction par suppression de $A_{0}$ peut se
faire par suppression de n'importe quel point $A_{p}$ en remarquant que si $%
(A_{n})_{n\in \Z}$ est un $CPD,$ alors $(A_{n+p})_{n\in \Z}$ en est un
\'{e}galement.

\section{Une propri\'{e}t\'{e} des $CPD$.}

Soit $\mathcal{A}=(A_{k})_{k\in \Z}$ un $CPD$ \`{a} $n$ sommets, avec $n>4$.
Tant que $n>4$, si deux c\^{o}t\'{e}s cons\'{e}cutifs sont dans le m\^{e}me
quadrant, on peut les remplacer, en appliquant l'op\'{e}ration de
r\'{e}duction ci-dessus, par un seul en supprimant leur point commun. Cet
algorithme se termine si $n=4$ et dans ce cas, le $CPD$ peut avoir deux
c\^{o}t\'{e}s cons\'{e}cutifs dans le m\^{e}me quadrant, ou si $n>4$ et le $%
CPD$ est tel que deux c\^{o}t\'{e}s cons\'{e}cutifs ne sont jamais dans le
m\^{e}me quadrant.

Nous allons d\'{e}montrer que cette derni\`{e}re situation est
contradictoire.

Soit donc $n$ un entier $>4$ et $\mathcal{A}=(A_{k})_{k\in \Z}$ un $CPD$
\`{a} $n$ sommets tel que, pour tout entier $k\in \Z,\;\overrightarrow{%
A_{k-1}A_{k}}$ et $\overrightarrow{A_{k}A_{k+1}}$ ne sont jamais dans le
m\^{e}me quadrant. Voici la premi\`{e}re cl\'{e} : Puisque le nombre de
c\^{o}t\'{e}s est sup\'{e}rieur ou \'{e}gal \`{a} cinq et qu'il n'y a que
quatre quadrants, il y a forc\'{e}ment deux vecteurs distincts $%
\overrightarrow{A_{p}A_{p+1}}$ et $\overrightarrow{A_{q}A_{q+1}}$ dans le
m\^{e}me quadrant qu'on peut supposer \^{e}tre le premier sans perte de
g\'{e}n\'{e}ralit\'{e}. On peut supposer, quite \`{a} renum\'{e}roter le $%
CPD,$ que $p=0$ et $1\le q\le n-1$. Mais la valeur $q=1$ est impossible
puisqu'alors $A_{0}A_{1}$ et $A_{q}A_{q+1}$ seraient deux c\^{o}t\'{e}s
cons\'{e}cutifs. De m\^{e}me $q=n-1 $ est impossible puisqu'alors $%
A_{q+1}=A_{n}=A_{0}$ et $A_{q}A_{q+1}$ et $A_{0}A_{1}$ seraient deux
c\^{o}t\'{e}s cons\'{e}cutifs. Ainsi l'entier $q$ v\'{e}rifie la condition $%
2\le q\le n-2$.

Le lemme suivant constitue la seconde cl\'{e}. Nous en avions fait une
d\'{e}monstration g\'{e}om\'{e}trique. La d\'{e}monstration analytique plus
simple que nous donnons nous a \'{e}t\'{e} inspir\'{e} par Fran\c{c}ois
Moulin.

\begin{Lem}
Avec les notations qui pr\'{e}c\'{e}dent, les vecteurs $\overrightarrow{%
A_{0}A_{1}}$ et $\overrightarrow{A_{1}A_{q}}.$ ne sont pas dans le m\^{e}me
quadrant.
\end{Lem}

\begin{Dem}
Nous supposons comme ci-dessus que $\overrightarrow{A_{0}A_{1}}$ est dans le
premier quadrant. Nous allons utiliser la m\'{e}thode analytique qui nous a
d\'{e}j\`{a} r\'{e}ussi. Soit donc $(x,y)$ les coordonn\'{e}es de $%
\overrightarrow{A_{0}A_{1}}$, celles de $\overrightarrow{A_{1}A_{2}}$
\'{e}tant $(x^{\prime },y^{\prime })$. On sait d\'{e}j\`{a} que ce dernier
vecteur n'est pas dans le premier quadrant, donc que l'une de ses deux
coordonn\'{e}es est strictement n\'{e}gative.\newline
Par hypoth\`{e}se $\det ($ $\overrightarrow{A_{0}A_{1}}$ ,$\overrightarrow{%
A_{1}A_{2}})=xy^{\prime }-x^{\prime }y>0$. Si $x^{\prime }\ge 0,$ alors $%
y^{\prime }<0$ et $xy^{\prime }-x^{\prime }y\le 0,$ ce qui est impossible,
donc $x^{\prime }<0$ et on sait d\'{e}j\`{a} que $x\ge 0$.\newline
Pla\c{c}ons-nous maintenant dans le rep\`{e}re direct $(A_{1},%
\overrightarrow{A_{1}A_{2}},\overrightarrow{A_{1}A_{0}})$. On sait
(d\'{e}monstration du th\'{e}or\`{e}me \ref{t1}) que les coordonn\'{e}es $%
(a,b)$ de $A_{q}$ dans ce rep\`{e}re sont strictement positives. Exprimons
le vecteur $\overrightarrow{A_{1}A_{q}}$ dans la base $(i,j)$. On a : 
\[
\overrightarrow{A_{1}A_{q}}=a\overrightarrow{A_{1}A_{2}}+b\overrightarrow{%
A_{1}A_{0}}=a(x^{\prime }i+y^{\prime }j)+b(-xi-yj)=(ax^{\prime
}-bx)i+(ay^{\prime }-by)j. 
\]
L'abscisse $ax^{\prime }-bx$ de ce vecteur est clairement strictement
n\'{e}gative. Ce vecteur ne peut donc \^{e}tre dans le premier quadrant, ce
qui termine la d\'{e}monstration du lemme.
\end{Dem}

En consid\'{e}rant les deux vecteurs $\overrightarrow{A_{q}A_{q+1}}$ et $%
\overrightarrow{A_{n}A_{n+1}}=\overrightarrow{A_{p}A_{p+1}}\ $et en
renum\'{e}rotant le $CPD$ \`{a} partir de $q$, on a \'{e}galement
d\'{e}montr\'{e} que le vecteur $\overrightarrow{A_{q+1}A_{p}}$ n'est pas
dans le premier quadrant.

Vu ce lemme, si $n>4,$ on peut faire des r\'{e}ductions successives par
suppression des sommets diff\'{e}rents de $A_{p},A_{p+1,}A_{q},A_{q+1}$ et
arriver au $CPD$ $(A_{p}A_{p+1}A_{q}A_{q+1})$ qui a quatre sommets et qui
v\'{e}rifie la condition : les deux vecteurs $\overrightarrow{A_{p}A_{p+1}}$
et $\overrightarrow{A_{q}A_{q+1}}$ sont dans le m\^{e}me quadrant tandis que
les deux autres vecteurs $\overrightarrow{A_{p+1}A_{q}}$ et $\overrightarrow{%
A_{q+1}A_{p}}$ n'y sont pas.

Le lemme suivant qui constitue la troisi\`{e}me cl\'{e} \'{e}tablit
l'impossibilit\'{e} d'une telle situation :\newline

\begin{Lem}
\label{l1}Soit $(A_{k})_{k\in \Z}$ un convexe polygonal direct \`{a} quatre
sommets tel que deux vecteurs cons\'{e}cutifs $\overrightarrow{A_{k-1}A_{k}}$
et $\overrightarrow{A_{k}A_{k+1}}$ ne sont jamais dans le m\^{e}me quadrant.
Alors les vecteurs $\overrightarrow{A_{0}A_{1}}$ (resp. $\overrightarrow{%
A_{1}A_{2}}$) et $\overrightarrow{A_{2}A_{3}}$ (resp. $\overrightarrow{%
A_{3}A_{4}}$) ne sont jamais dans le m\^{e}me quadrant.
\end{Lem}

\begin{Dem}
Notons $(A,B,C,D)$ les quatres sommets $(A_{0},A_{1},A_{2},A_{3})$ du
quadrilat\`{e}re et raisonnons par l'absurde en supposant que les vecteurs $%
AB=(x,y)$ et $CD=(x^{\prime \prime },y^{\prime \prime })$ sont dans le 1er
quadrant (i.e. $x,x^{\prime \prime },y,y^{\prime \prime }\ge 0$). notons $%
BC=(x^{\prime },y^{\prime })$ et cherchons \`{a} le localiser. On sait
d\'{e}j\`{a} qu'il ne peut \^{e}tre dans le 1er ou le 4eme quadrant. Il ne
peut \^{e}tre dans le 2eme quadrant (i.e. $x^{\prime }\le 0,\;y^{\prime }\ge
0$) car dans ce cas on aurait : 
\[
\det (BC,CD)=x^{\prime }y^{\prime \prime }-x^{\prime \prime }y^{\prime }\le
0. 
\]
On a donc $y^{\prime }<0$ et $BC$ est strictement dans le 3eme quadrant ($%
x^{\prime }<0$ et $y^{\prime }<0$).\newline
Le m\^{e}me raisonnement appliqu\'{e} au vecteur $CD$ montre que $DA$ est
aussi strictement dans le 3eme quadrant. $CB$ et $AD$ sont ainsi dans le 1er
quadrant.\newline
On a ainsi \'{e}tabli que les vecteurs $AB,AD,CB,CD$ sont dans le 1er
quadrant. Les points $B$ et $D$ sont donc dans le 1er quadrant d'origine $C$
et aussi dans le 1er quadrant d'origine $A$. Dans ces conditions, le point $%
A $ ne peut \^{e}tre dans le 2eme 3eme ou 4eme quadrant d'origine $C$ sinon
les deux segments $[A,C]$ et $[B,D]$ ne peuvent se couper qu'aux points $%
C,\;B$, $D\ $ou $A,$ ce qu'interdit le lemme des diagonales. Le point $A$
est donc dans le 1er quadrant d'origine $C$. Mais alors, $C$ est dans le
3eme quadrant d'origine $A$ et on retombe sur le m\^{e}me interdit. Dans
tous les cas, l'appartenance de $AB$ et $CD$ au m\^{e}me quadrant est
contradictoire.
\end{Dem}

Nous avons ainsi d\'{e}montr\'{e} le th\'{e}or\`{e}me suivant :

\begin{Th}
Soit $(A_{k})_{k\in \Z}$ un convexe polygonal direct \`{a} $n$ sommets. On
suppose que pour tout $k\in \Z,$ les vecteurs $\overrightarrow{A_{k-1}A_{k}}$
et $\overrightarrow{A_{k}A_{k+1}}$ ne sont jamais dans le m\^{e}me quadrant.
Alors $n\le 4$.
\end{Th}

Dans ce qui suit, nous allons appliquer ce th\'{e}or\`{e}me \`{a} un
ensemble de polygones convexes produits \`{a} partir d'une famille finie de
droites que nous allons d\'{e}finir.

\section{D\'{e}finition et \'{e}tude de $D_{m,n}$ et de $CF(m,n)$}

Soit $(m,n)\in \N^{*2}.$ On consid\`{e}re l'ensemble $D_{m,n}$ des droites
rencontrant le carr\'{e} unit\'{e} $CU=[0,1]^{2}$ et d'\'{e}quation$\;\fbox{$%
ux+vy-w=0$}$\ dont les coefficients $u,v,w$ sont des \underline{entiers
relatifs} et v\'{e}rifient les conditions :$\;\;$%
\[
-m\le u\le m,\;-n\le v\le n,\;(u,v)\ne (0,0),\;w\in \Z
. 
\]
Les droites d'\'{e}quation $x=0,x=1,y=0$ et $y=1$ appartiennent bien
\'{e}videmment \`{a} $D_{m,n},$ ainsi que les droites d'\'{e}quation $-x+y=0$
et $x+y-1=0.$ Notons que l'ensemble $D_{m,n}$ ainsi d\'{e}fini est
constitu\'{e} d'un nombre fini de droites. En effet, si la droite $D$
d'\'{e}quation $ux+vy-w=0$ est dans $D_{m,n}$ elle contient un point $M$ de
coordonn\'{e}es $(a,b)\in CU.$ Ceci montre que l'entier $w$ satisfait
l'in\'{e}galit\'{e} : 
\[
\left| w\right| \le \left| u\right| a+\left| v\right| b\le m+n. 
\]
Pour un couple $(u,v)$ donn\'{e}, il n'y a donc qu'un nombre fini de valeurs
de $w$ qui conviennent. Il suffit alors de remarquer que les couples $(u,v)$
admissibles sont en nombre fini. On obtient accessoirement une majoration du
cardinal de $D_{m,n}$ par $(2m+1)(2n+1)(2m+2n+1).$ On peut s'int\'{e}resser
au cardinal de $D_{m,n}$ et le calculer, mais c'est une autre histoire.

\begin{Def}[Complexe de Farey]
On appelle complexe de Farey d'ordre $(m,n)$ le compl\'{e}mentaire de $%
\bigcup D_{m,n}\ $dans le carr\'{e} unit\'{e} $CU.$ On le note $CF(m,n);$
\end{Def}

L'ensemble $CF(m,n)\ $est un ouvert du plan. Soit $K$ l'une de ses
composantes connexes. Par d\'{e}finition elle est non vide born\'{e}e et
aucune des droites de $D_{m,n}$ ne peut rencontrer $K$, si bien que $K$ est
toujours contenu dans l'un des deux demi-plans ouverts d\'{e}termin\'{e} par
chacune des droites de $D_{m,n}.$ C'est ainsi que $K$ est contenu dans
l'intersection $H$ d'une famille finie de demi-plans ouverts qui sont, on le
sait, convexes. $H$ est ainsi une partie convexe de $CF(m,n)$, donc connexe.
On a donc, par maximalit\'{e} de $K$, $K=H$. L'adh\'{e}rence de $K$ est donc
un polygone convexe. Il existe donc un $CPD$ $\mathcal{A}=(A_{k})_{k\in \Z}$
\`{a} $p$ c\^{o}t\'{e}s tel que $LP(\mathcal{A})=Fr(K)$. On rappelle que $LP(%
\mathcal{A})$ est la r\'{e}union des segments $[A_{k},A_{k+1}]$ et que $%
Fr(K) $ est la fronti\`{e}re de $K$.

Lorsqu'on fait une simulation graphique en dessinant les $D_{m,n}$ pour les
petites valeurs de $m$ et $n,$ (voir la figure en fin de texte), on s'aper\c{%
c}oit que les $Fr(K)$ sont des triangles ou des quadrilat\`{e}res. Le
calcul, pour ces m\^{e}mes petites valeurs, le confirme. Nous nous proposons
de d\'{e}montrer que $p\in \{3,4\},$ c'est-\`{a}-dire que $Fr(K)$ est
toujours un triangle ou un quadrilat\`{e}re. Cela nous donnera l'occasion
d'en d\'{e}montrer un peu plus sur ces triangles et quadrilat\`{e}res.

\subsection{Une propri\'{e}t\'{e} de $D_{m,n}$}

Notons d'abord que l'ensemble $D_{m,n}$ est invariant par les deux
sym\'{e}tries $\sigma _{x}$ et $\sigma _{y}$ par rapport aux droites $y=1/2$
et $x=1/2$. En effet si $D\in D_{m,n}$ a pour \'{e}quation $ux+vy-w=0$, $%
\sigma _{x}(D)$ a pour \'{e}quation $ux+v(1-y)-w=ux-vy+v-w=0$ et $\sigma
_{y}(D)$ a pour \'{e}quation $u(1-x)+vy-w=-ux+vy+u-w=0,$ et ces deux droites
appartiennent \`{a} $D_{m,n}$.

Soit maintenant $A,B$ deux points distincts de $CU$ tels que le vecteur $%
AB=(a,b)$ soit dans le 2eme quadrant $(a\le 0,\;b\ge 0)$. Alors le vecteur $%
\sigma _{y}(A)\sigma _{y}(B)$ de coordonn\'{e}es $(-a,b)$ appartient au 1er
quadrant. De m\^{e}me si ce vecteur appartient au 4eme quadrant $(a\ge
0,b\le 0),$ le vecteur $\sigma _{x}(A)\sigma _{x}(B)$ de coordonn\'{e}es $%
(a,-b)$ appartient au 1er quadrant. Enfin, s'il appartient au 3eme quadrant $%
(a\le 0,b\le 0)$, le vecteur $\sigma _{y}\circ \sigma _{x}(A)\sigma
_{y}\circ \sigma _{x}(B)$ de coordonn\'{e}es $(-a,-b)$ appartient au 1er
quadrant. Ces invariants vont nous permettre une simplification dans la
d\'{e}monstration du lemme suivant qui constitue la quatri\`{e}me cl\'{e} de
notre d\'{e}monstration :

\begin{Lem}[des trois points]
Soient $A,B,C$ trois points distincts non align\'{e}s du carr\'{e} unit\'{e} 
$CU.$ On suppose que\ :

\begin{enumerate}
\item  {Les vecteurs $\overrightarrow{AB}$ et $\overrightarrow{BC}$ sont
dans le m\^{e}me quadrant. }

\item  {Les droites $(A,B)$ et $(B,C)$ sont dans $D_{m,n}$ }
\end{enumerate}

Alors il existe une droite $D\in D_{m,n}$ d'\'{e}quation $\phi (M)=0$ qui
passe par le point $B$ et qui traverse le segment $AC$ (i.e. $\phi (A)\phi
(C)<0$).
\end{Lem}

\begin{Dem}
Nous supposerons que les vecteurs $\overrightarrow{AB}$ et $\overrightarrow{%
BC}$ sont dans le premier quadrant, les autres cas s'y ramenant par les
sym\'{e}tries invoqu\'{e}es ci-dessus. Nous allons trouver $\phi $ comme
diff\'{e}rence d'une \'{e}quation de $(A,B)$ et d'une \'{e}quation de $(B,C)$
bien choisies.

Notons $(x_{M},y_{M})$ les coordonn\'{e}es d'un point $M$ du plan. Les
coordonn\'{e}es des points $A,B$ et $C,$ vue la premi\`{e}re hypoth\`{e}se,
v\'{e}rifient les in\'{e}galit\'{e}s : 
\[
x_{A}\le x_{B}\le x_{C}\;\;\;\text{et}\;\;\;y_{A}\le y_{B}\le y_{C} 
\]

Les droites $(A,B)$ et $(B,C)$ admettent respectivement des \'{e}quations de
la forme $f(M)=ux+vy-w=0$ et $f^{\prime }(M)=u^{\prime }x+v^{\prime
}y-w^{\prime }=0\ $o\`{u} $M=(x,y)$ est un point g\'{e}n\'{e}rique du plan.
Ces droites \'{e}tant dans $D_{m,n},$ leurs coefficients peuvent \^{e}tre
choisies comme des entiers v\'{e}rifiant les in\'{e}galit\'{e}s : 
\[
|u|\le m,\;|u^{\prime }|\le m,\;|v|\le n,\;|v^{\prime }|\le n.\; 
\]
Les entiers $u$ et $v$ ne sont jamais simultan\'{e}ment nuls. Il en est de
m\^{e}me des entiers $u^{\prime }$ et $v^{\prime }$.

La droite $(A,B)$ admet \'{e}galement pour \'{e}quation : 
\begin{eqnarray*}
g(M) &=&\det (AB,AM) \\
&=&(x_{B}-x_{A})(y-y_{A})-(y_{B}-y_{A})(x-x_{A}) \\
&=&\alpha x+\beta y-\gamma =0.
\end{eqnarray*}
On a $\beta =(x_{B}-x_{A})\ge 0$ et $\alpha =-(y_{B}-y_{A})\le 0$.\ De plus
ces deux r\'{e}els $\alpha $ et $\beta $ ne sont pas nuls simultan\'{e}ment
car $A\ne B$. $f$ et $g,$ \'{e}quations de la m\^{e}me droite $(A,B)$ sont
donc proportionnels, avec un coefficient de proportionnalit\'{e} non nul. En
changeant \'{e}ventuellement $f$ en $-f$ on peut supposer que ce coefficient
est strictement positif. Ainsi on aura $u\le 0$ et $v\ge 0$.

On peut faire de m\^{e}me avec $f^{\prime }$ en posant $g^{\prime }(M)=\det
(BC,BM)$ et faire en sorte que $u^{\prime }\le 0$ et $v^{\prime }\ge 0$.

Etudions maintenant le signe de $g(C)$ et celui de $g^{\prime }(A)$ qui ne
sont jamais nuls puisque, par hypoth\`{e}se $A,B,C$ ne sont pas align\'{e}s.
On a : 
\begin{eqnarray*}
g(C) &=&\det (AB,AC)=\det (AB,AB+BC)=\det (AB,BC) \\
g^{\prime }(A) &=&\det (BC,BA)=-\det (BA,BC)=\det (-BA,BC)=\det (AB,BC).
\end{eqnarray*}
Ces d\'{e}terminants (\'{e}gaux et non nuls) sont donc de m\^{e}me signe.
Ainsi $f(C)$ et $f^{\prime }(A)$ sont de m\^{e}me signe (et non nuls tous
les deux). Posons : 
\[
\phi (M)=f(M)-f^{\prime }(M)=(u-u^{\prime })x+(v-v^{\prime })y-(w-w^{\prime
}). 
\]
Vu les signes de $u,u^{\prime }$ et ceux de $v,v^{\prime },$ et la
stabilit\'{e} de $\Z$ par soustraction, $\phi $ est l'\'{e}quation d'une
droite de $D_{m,n}$ qui passe par le point $B$. On a de plus : 
\[
\phi (A)=-f^{\prime }(A)\;\;\;\text{et\ \ \ }\phi (C)=f(C)\;\;\;\text{donc\
\ \ }\phi (A)\phi (C)<0.\ 
\]
La droite $D$ d'\'{e}quation $\phi (M)=0$ est une solution au probl\`{e}me
pos\'{e} CQFD.
\end{Dem}

La conclusion de ce lemme peut se formuler autrement : le point $A$ est dans
l'un des demi-plans ouverts d\'{e}termin\'{e}s par la droite $D$ et le point 
$C$ est dans l'autre. Ce lemme nous permet de donner le th\'{e}or\`{e}me
suivant qui est la propri\'{e}t\'{e} essentielle de $Fr(K)$.

\begin{Th}
\label{T2}Soit $K$ une composante connexe de $CF(m,n)$ et soit $A,B,C$ trois
sommets cons\'{e}cutifs de $Fr(K)$ (qui est un convexe polygonal direct).
Alors les vecteurs $\overrightarrow{AB}$ et $\overrightarrow{BC}$ ne sont
jamais dans le m\^{e}me quadrant.
\end{Th}

\begin{Dem}
Raisonnons par l'absurde et supposons que $AB$ et $BC$ sont dans le m\^{e}me
quadrant. Alors d'apr\`{e}s le lemme pr\'{e}c\'{e}dent, il existe une droite 
$D\in D_{m,n}$ passant par $B,$ d'\'{e}quation $\phi ,$ et telle que $\phi
(A)\phi (C)<0.$ On peut, en changeant \'{e}ventuellement $\phi $ en $-\phi ,$
supposer que $\phi (A)>0$ et $\phi (C)<0.$ les points $A$ et $C$ \'{e}tant
sur la fronti\`{e}re de $K$ et $\phi $ \'{e}tant continue, il existe un
point $A^{\prime }$ de $K$ tel que $\phi (A^{\prime })>\phi (A)/2$ et un
point $C^{\prime }$ de $K$ tel que $\phi (C^{\prime })<\phi (C)/2$.\newline
Notons $P_{+}$ le demi-plan ouvert $\{M\;|\;\phi (M)>0\}$ et $P_{-}$ le
demi-plan ouvert$\{M\;|\;\phi (M)<0\}.$ Comme , par d\'{e}finition de $K,$ $%
D\cap K=\emptyset $, on a 
\[
K=(P_{+}\cap K)\cup (P_{-}\cap K) 
\]
Ces deux ouverts sont non vides (le premier contient $A^{\prime }$, le
second contient $C^{\prime }$) et disjoints. $K$ est donc non connexe, ce
qui est une contradiction.\newline
\end{Dem}

\subsection{Les composantes connexes de $CF(m,n)$.}

Rappelons que $CF(m,n)$, contenu dans le carr\'{e} unit\'{e}, est born\'{e}.
Soit $K$ l'une de ses composantes connexes (qui est donc born\'{e}e). Le
th\'{e}or\`{e}me pr\'{e}c\'{e}dent nous apprend que $Fr(K)$ est un $LP(%
\mathcal{A})$ d'un $LPC$ $\mathcal{A}$ qui v\'{e}rifie la propri\'{e}t\'{e}
: deux c\^{o}t\'{e}s cons\'{e}cutifs ne sont pas dans le m\^{e}me quadrant.
Nous avons d\'{e}montr\'{e} plus haut qu'un tel $LPC$ ne pouvait avoir plus
de quatre c\^{o}t\'{e}s. Comme l'adh\'{e}rence de $K$ est un polygone
convexe, $\mathcal{A}$ a trois c\^{o}t\'{e}s au moins. Nous avons ainsi
d\'{e}montr\'{e} le th\'{e}or\`{e}me suivant :

\begin{Th}[Tajine-Daurat]
L'adh\'{e}rence de toute composante connexe du complexe de Farey $CF(m,n)$
est un triangle ou un quadrilat\`{e}re.
\end{Th}

Nous pouvons am\'{e}liorer le r\'{e}sultat de ce th\'{e}or\`{e}me en
pr\'{e}cisant la nature de $Fr(K)$ selon la position des c\^{o}t\'{e}s
cons\'{e}cutifs. On est dans l'une des deux configurations suivantes.

- Si $Fr(K),$ en tant que $CPD,$ est tel que deux c\^{o}t\'{e}s
cons\'{e}cutifs sont dans deux quadrants oppos\'{e}s (i.e. $Q_{1}$ et $Q_{3}$
ou $Q_{2}$ ou $Q_{4}$), alors $Fr(K)$ est un triangle.

- Sinon $Fr(K)$ est un quadrilat\`{e}re. Si on en prend un repr\'{e}sentant
convexe polygonal direct $(A_{n})_{n\in \Z}$ alors il existe un entier $p$
tel que, pour $i\in [1..4],\;\overrightarrow{A_{p+i-1}A_{p+i}}$ est dans le
quadrant $Q_{i}$.\newline

La d\'{e}monstration de ces deux r\'{e}sultats est analogue \`{a} celle du
lemme \ref{l1}. Elle repose sur le fait simple suivant :

\begin{Lem}
Si $Fr(K)$ est un quadrilat\`{e}re repr\'{e}sent\'{e} par le convexe
polygonal direct $ABCD,$ alors la condition ``$AB\in Q_{1}$ et $BC\in Q_{3}$%
'' est contradictoire.
\end{Lem}

\begin{Dem}
En effet $CD$ ne peut appartenir \`{a} $Q_{2}$ pour une raison
d'orientation. Par le lemme \ref{l1}, il ne peut appartenir \`{a} $Q_{1}$.
Par le th\'{e}or\`{e}me \ref{T2}, il ne peut appartenir \`{a} $Q_{3}$. Il
appartient donc \`{a} $Q_{4}$. Pour des raisons analogues $DA$ appartient
\`{a} $Q_{2}$.\newline
Ceci montre que $B$ et $D$ sont dans le demi-plan $x\ge 0$ d'origine $A$ et
aussi dans le demi-plan $x\ge 0$ d'origine $C$, ce qui d\'{e}montre que les
diagonales $AC$ et $BD$ ne peuvent se traverser, contredisant le lemme des
diagonales.
\end{Dem}

\section{Conjecture forte de Tajine-Daurat et remerciements.}

Nous conjecturons que les composantes connexes born\'{e}es du
compl\'{e}mentaire de $\bigcup D_{m,n}$ dans le plan sont des triangles ou
des quadrilat\`{e}res. Cette propri\'{e}t\'{e} peut-elle \^{e}tre
\'{e}tendue aux composantes connexes non born\'{e}es et quel sens cette
extension aurait-elle ?

Pour conclure, nous tenons \`{a} remercier Fran\c{c}ois Rideau et Fran\c{c}%
ois Moulin pour leurs contributions, Anatole Kh\'elif qui a manifest\'{e} un
grand int\'{e}r\^{e}t pour ce travail, ainsi que Mohamed Tajine qui a fait
un travail de relecture tr\`{e}s attentif et qui nous a fourni le source du
dessin incorpor\'{e} \`{a} ce texte.

{\Large \hskip-1cm \textbf{R\'{e}f\'{e}rences :} }

{[1] Marcel Berger. G\'{e}om\'{e}trie livre 3: Convexes et polytopes,
poly\`{e}dres r\'{e}guliers, aires et volumes. CEDIC/Fernand Nathan, 1977. }

{[2] Malcolm Douglas McIlroy. A note on discrete representation of lines.
ATT Technical Journal, 64(2) :481-490, 2984. }

{[3] Alain Daurat, Mohamed Tajine, Mahdi Zouaoui. About the frequencies of
some patterns in digital planes. application to area estimators. Computer
Graphics, 2008. }\newline

\newpage%

\FRAME{itbpFU}{1.9666in}{1.9934in}{0in}{\Qcb{Cas $m=4,\;n=3$}}{}{%
diagramme_farey.ps}{\special{language "Scientific Word";type
"GRAPHIC";maintain-aspect-ratio TRUE;display "USEDEF";valid_file "F";width
1.9666in;height 1.9934in;depth 0in;original-width 228.75pt;original-height
231.8125pt;cropleft "0";croptop "1";cropright "1";cropbottom "0";filename
'C: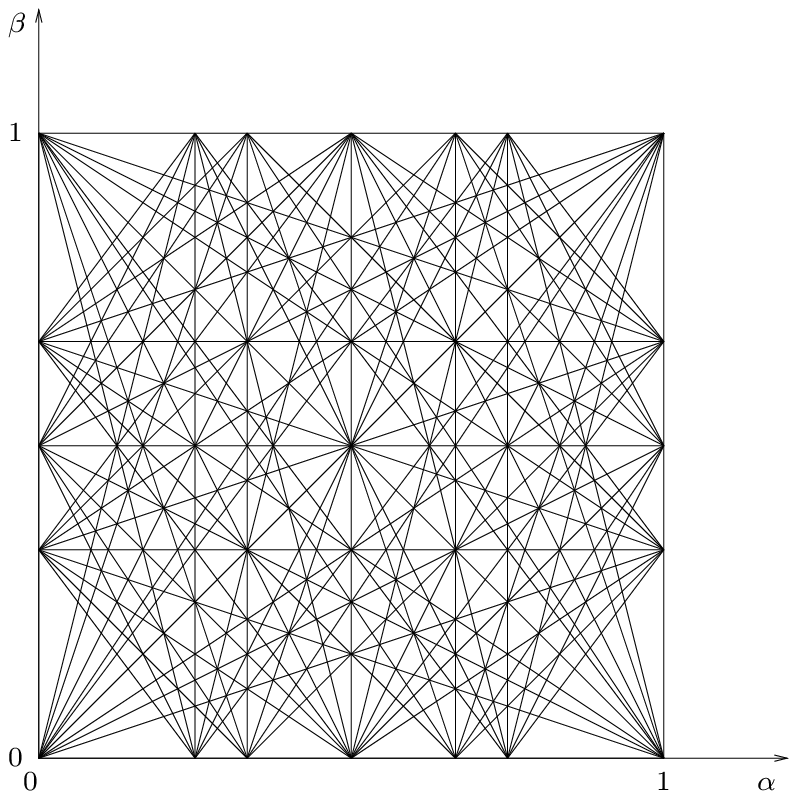';file-properties "XNPEU";}}

\end{document}